\documentclass[letterpaper, 10 pt, conference]{cssconf}  %
\IEEEoverridecommandlockouts     
\usepackage{multirow}
\usepackage{scrextend}
\usepackage{fancyhdr}
\usepackage{balance}
\usepackage{mathtools}
\usepackage{comment}

\usepackage{tikz}
\usepackage{pgfplots}
\usepackage[compress]{cite}
\usepackage{amsmath}

\usepackage[mathscr]{euscript}
 \let\mathscr\relax% just so we can load this and rsfs
\usepackage[scr]{rsfso}
\usepackage{amsthm}
\usepackage{amssymb}
\usepackage{multicol}
\usepackage[colorlinks=true, pdfstartview=FitV, linkcolor=blue,
citecolor=blue, urlcolor=blue]{hyperref}
\usepackage{amsmath}

\newcommand{\R}{\mathbb{R}}

\newcommand{\bm}{\mathbb{}}

\newtheorem{lemma}{Lemma}
\newtheorem{remark}{Remark}
\newtheorem{proposition}{Proposition}
\title{\LARGE \bf
Data-Driven Adjustable Robust Optimization 
}

 \author{
 	 Xiaoxing Ren,  Alessio Moreschini, Zhongda Chu, Yulong Gao, Thomas Parisini
   \thanks{
		This work has been partially supported by the European Union’s Horizon 2020 Research and Innovation Programme under grant agreement No. 739551 (KIOS CoE).
	}
 	\thanks{X. Ren, A. Moreschini, Z. Chu, Y. Gao, and T. Parisini are with the Department of Electrical and Electronic Engineering, Imperial College London, UK. 
  T. Parisini is also with the Department of Electronic Systems, Aalborg University, Denmark, and with the Department of Engineering and Architecture, University of Trieste, Italy
  {\tt\small \{xiaoxing.ren, a.moreschini, z.chu18, yulong.gao,t.parisini\}@imperial.ac.uk}}
 	}

\begin{document}

\maketitle

% don't print page numbers
\thispagestyle{empty}
\pagestyle{empty}
% EVERYTHING ABOVE THIS LINE IS JUST PREABLE, NO NEED TO MESS WITH IT.__________________________________________________________________________________________
%
\begin{abstract}
In this paper, we develop a two-stage
data-driven approach to address the adjustable robust optimization problem, where the uncertainty set is adjustable to manage infeasibility caused by significant or poorly quantified uncertainties. 
In the first stage, we synthesize an uncertainty set to ensure the feasibility of the problem as much as possible using the collected uncertainty samples. In the second stage, we find the optimal solution while ensuring that the constraints are satisfied under the new uncertainty set. 
This approach enlarges the feasible state set, at the expense of the risk of possible constraint violation. 
We analyze two scenarios: one where the uncertainty is non-stochastic, and another where the uncertainty is stochastic but with unknown probability distribution, leading to a distributionally robust optimization problem. 
In the first case, we scale the uncertainty set and find the best subset that fits the uncertainty samples. 
%reformulate the problem for several typical uncertainty sets. 
In the second case, we employ the Wasserstein metric to quantify uncertainty based on training data, and for polytope uncertainty sets, we further provide a finite program reformulation of the problem. %\yulong{Do we provide a tractable reformulation?}
The effectiveness of the proposed methods is demonstrated through an optimal power flow problem.
\end{abstract}

\section{Introduction} \label{se:introduction}
Many real-world decision problems arising in engineering and management involve uncertain parameters,
where the uncertainty model is not stochastic but rather deterministic and set-based. 
Robust optimization is a popular approach to optimization with set-based uncertainty. 
In contrast to stochastic optimization, which aims to find a decision that minimizes the expected cost based on the probability distribution of uncertain parameters,  in robust optimization, the decision-maker constructs a solution that is feasible for any realization of the uncertainty in a given set. 

The general robust optimization problem takes on the form 
\begin{equation} \label{eq:problem_formulation}
\begin{array}{ll}
\underset{x}{\min } & g(x)  
\\
\text {s.t.} 
& f(x, \xi) \leq 0, \quad x \in \mathcal{X},	\quad \forall \xi \in \mathbb{S},
\end{array}
\end{equation}
 where $x \in \mathcal{X}$, $\mathcal{X} \subset \mathbb{R}^n$ is the decision variable, and $\xi \in \R^m$ is the uncertain parameter, which is known to belong to a set $\mathbb{S}$ (called the {\em uncertainty set}). The goal of Problem~\eqref{eq:problem_formulation} is to compute the minimum cost solution among all those solutions which are feasible for all realizations of the uncertainty within the set $ \mathbb{S}$.
It has been widely applied in various fields, such as portfolio optimization, learning, control, estimation and power systems~\cite{beyer2007robust,bertsimas2011theory,saltik2018outlook,conejo2022robust,ning2019optimization}. 

In the past decade, the very significant availability of data has driven the development of data-centered approaches in robust optimization.
%, in particular, when the uncertainty set is unknown, design data-driven uncertainty sets.
When the uncertainty is stochastic and probability distribution is unknown, the approach is referred to as {\em distributionally robust optimization}~\cite{van2021data,jiang2024distributionally,levy2020large,lin2022distributionally}, where the optimal solutions are evaluated under the worst-case expectation with respect to a family of probability distributions of the uncertain parameters, which are called ambiguity sets.
Data-driven and learning-based methods are crucial for leveraging finite training data to define uncertainty or ambiguity sets. This approach ensures robust decision-making on unseen samples, thereby enhancing the reliability and robustness of optimization solutions under distributional ambiguity, see, e.g.,~\cite{margellos2014road,bertsimas2018data,kuhn2019wasserstein,alexeenko2020nonparametric}.
% On the other hand, when the uncertainty is non-stochastic,~\cite{bertsimas2018data} designs uncertainty sets from data using hypothesis tests,
%~\cite{alexeenko2020nonparametric}
In~\cite{margellos2014road}, a data-driven approach is developed to estimate intervals that cover each component of a random vector with high probability, resulting in hyper-rectangular uncertainty. In~\cite{bertsimas2018data}, uncertainty sets are designed from data using statistical hypothesis tests. The data-driven approach in ~\cite{alexeenko2020nonparametric} constructs uncertainty sets with the guarantee that their probability mass is within a given tolerance of the target mass with high confidence.
Besides, ~\cite{kuhn2019wasserstein} proposes data-driven methods to construct the ambiguity sets using Wassertein metrics and seek decisions that perform best in view of the worst-case distribution within the ambiguity sets.
In traditional robust optimization, the worst-case viewpoint of uncertainty is considered, which leads to overly conservative solutions and the problem may well become infeasible. 

\smallskip

\textbf{Motivating Example:} Consider optimization problem~\eqref{eq:problem_formulation} with 
%constraint function 
$f(x, \xi) = x + \xi - 5$, 
%optimization set 
$\mathcal{X}=[1,2]$, 
%and 
%uncertainty set 
$\mathbb{S} = \{ \xi \mid  3 \leq \xi \leq 5 \}$.
\begin{figure}[b!]
    \centering
    \vspace{-1em}
    \includegraphics[width=0.475\linewidth]{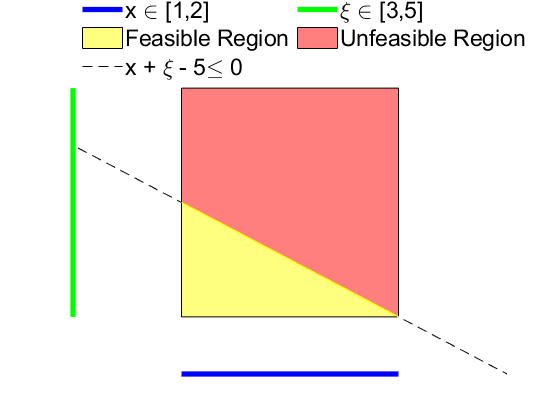}
    \includegraphics[width=0.475\linewidth]{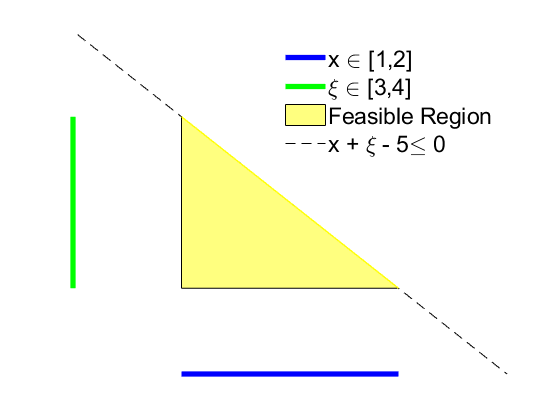}
    \vspace{-1em}
    \caption{Illustration of the feasible and unfeasible regions of the problem addressed in the motivating example. On the left, we depict the problem with its original uncertainty set, while on the right, we show the problem with a reduced uncertainty set. The benefit of the right-hand side problem is the elimination of the unfeasible~region.}
    \label{fig:motivating_example}
  %  \vspace{-1em}
\end{figure}
As shown in Fig.~\ref{fig:motivating_example}, 
it can be readily seen that for all $\xi\in(4,5]$ there is no feasible solution $x \in \mathcal{X}=[1,2]$ that allows the constraint $ f(x, \xi) \leq 0$ to hold. Nevertheless, to ensure the feasibility of the optimization problem at hand, one might consider a restriction to a subset of the uncertainty set $\hat{\mathbb{S}}\subset\mathbb{S}$, with $\hat{\mathbb{S}} = \{ \xi \mid  3 \leq \xi \leq 4 \}$. 
%See Fig.~\ref{fig:motivating_example}. 
\hfill $\blacktriangleleft$

\smallskip

Our goal is to solve the data-driven adjustable robust optimization problem with a priority on feasibility.
% by constructing uncertainty sets based on data samples drawn from an unknown distribution
We construct uncertainty sets based on uncertainty samples such that they can ensure the existence of a feasible solution under uncertainty.
In particular, inspired by~\cite{jaillet2022satisficing} we use a measure $\mu(\cdot)$ to choose the new uncertainty set and solve the following constrained optimization problem 
\begin{equation} \label{eq:penalty}
\begin{array}{ll}
\underset{x}{\max } & \underset{\mathbb{\hat{S}} \subset \mathbb{S}}{\max } \quad \mu(\mathbb{\hat{S}}) + \varrho(\hat{\Xi}\subseteq\mathbb{\hat{S}})
\\
\text {s.t.} 
& f(x, \xi) \leq 0, \ x \in \mathcal{X}, \ \quad \forall \xi \in \mathbb{\hat{S}} \subseteq \mathbb{S},
\end{array}
\end{equation}
where $\mu: \mathbb{\hat{S}} \rightarrow \mathbb{R} $ is a function defined on the domain $\mathbb{S} \subset \mathbb{R}^n$ ensuring that, if $\mu(\mathbb{\hat{S}}_1) > \mu(\mathbb{\hat{S}}_2)$, then $\mathbb{\hat{S}}_1$ is preferable than~$\mathbb{\hat{S}}_2$, whereas $\varrho: \mathbb{\hat{S}}  \rightarrow \mathbb{R}$ is a penalty function defined on the domain of the uncertainty samples set $\hat{\Xi}$.

When the uncertainty $\xi$ is non-stochastic, $\mu(\mathbb{{S}} 
) $ could be the size  or the volume of the adjustable uncertainty set, and we find a subset $\mathbb{\hat{S}} \in \mathbb{S}$ such that the  problem is feasible while maximizing $\mu(\mathbb{\hat{S}})$. On the other hand, when $\xi$ is stochastic with distribution $\mathbb{P}$, we find the set $\mathbb{\hat{S}}$ that maximizes $\mathbb{P}(\xi \in \mathbb{\hat{S}})$, i.e., the probability that an uncertain parameter  $\xi$  lies within the set $ \mathbb{\hat{S}}$. In practice, the distribution $\mathbb{P}$ is often unknown and can be approximated by observed uncertainty samples. This give rise to the adjustable distributional robust optimisation problem.
In particular, when the distribution  can be partially observed through a finite set of independent samples, we employ the uncertainty quantification framework proposed in~\cite{bertsimas2018data}. Then $\mu(\mathbb{\hat{S}})$ is related to $\mathbb{P}(\xi \in \mathbb{\hat{S}})$.

The main contributions are summarized as follows.

\emph{1)} We propose a two stage data-driven robust optimization algorithm. 
Different from standard robust optimization, where the uncertainty set is fixed, the proposed approach constructs a new smaller uncertainty set while ensuring feasibility. 
That is, the proposed algorithm gives priority to feasibility and achieves feasibility in an uncertain environment as well as possible. 

\emph{2)} When the uncertainty is non-stochastic,
we study uncertainty sets of practical interest, including norm-balls, ellipsoids and polyhedral sets, and demonstrate that these uncertainty sets can be easily adjusted.  
Unlike related works on adjustable uncertainty~\cite{zhang2017robust,gao2022robust,liu2023opportunistic,jaillet2022satisficing}, which compute the largest set according to a given metric that the system can tolerate without violating the constraints, we construct the largest uncertainty from the uncertainty samples.    
While ~\cite{gao2024learning} also approximates the uncertainty set using the collected disturbance samples, it does not address distributional uncertainty, which is essential for constructing a distributionally robust optimization framework. 
In this case, we utilize the Wasserstein metric to quantify uncertainty
based on training data. Additionally, when the uncertainty set is polytopic, the problem can be reformulated as a finite program.
% , providing rigorous guarantees for out-of-sample performance.
%

\emph{3)} We study the optimal power flow problem that arises in power systems, and show that it can be formulated as a robust optimal problem with an adjustable uncertainty set. 
We investigate cases where the uncertainty is non-stochastic and where the uncertainty is stochastic, demonstrating the performance guarantees of the proposed robust optimization scheme with adjustable uncertainty.

\textbf{Notations:}
The identity matrix is denoted by $\mathrm{I}$.
The vector with all elements equal to $1$ is denoted by $\mathrm{1}$.
A positive (semi-) definite matrix $Q$ is denoted by $Q \succ 0$ ($Q \succeq 0$). 
The cardinality of a set $\mathbb{X}$ is denoted by $|\mathbb{X}|$. 
The inner product of two vectors $\bm{a}, \bm{b} \in \mathbb{R}^n$ is denoted by $\langle \bm{a}, \bm{b} \rangle= \bm{a}^\top \bm{b}$. Given a norm $\|\cdot\|$ on $\mathbb{R}^m$, the dual norm is defined through $\| \bm{z} \|_* = \sup_{\|\bm{\xi}\| \leq 1} \langle \bm{z}, \bm{\xi} \rangle$, and the dual norm of $2$-norm is still $2$-norm.
Finally, the indicator function is denoted by $\mathbb{I}$, that is for a given set $\Xi$ one has that $\mathbb{I}(\xi \in \Xi) = 1$ if $\xi \in \Xi$, and  $\mathbb{I}(\xi \in \Xi) = 0$ otherwise. 

\section{Data-driven adjustable robust optimization}\label{sec:uniform_uncertain}
In this section, we consider Problem~\eqref{eq:penalty} where the uncertainty is non-stochastic. 
We construct a two-stage approach, the first stage is to find the largest set inside ${\mathbb{S}}$ that ensure feasibility with the given uncertainty samples. Then second stage is to achieve the best possible performance with the  uncertainty set obtained in the first stage.

Given the set of i.i.d. uncertainty samples $\hat{\Xi} =  \lbrace\hat{\xi}_i \rbrace_{i=1,...,N}$, the penalty term introduced in \eqref{eq:penalty} can be specified by the indicator function on the uncertainty samples set, that is $\varrho(\hat{\Xi}\subseteq\mathbb{\hat{S}}) \coloneqq  \gamma\sum_{i=1}^N \mathbb{I}(\hat{\xi}_i \in \hat{\mathbb{S}})$ for some given $\gamma > 0$.
This indicator function measures how well the new uncertainty set  $\mathbb{\hat{S}}$
  fits the uncertainty samples, specifically, we aim to include as many samples as possible. 
% \yulong{Can we add some justification of this item?}
With this in mind, the first stage is formulated as follows:

\begin{equation} \label{eq:stage_1}
\begin{array}{lll}
\text{\emph{Stage 1}:} & \underset{x}{\max } & \underset{\mathbb{\hat{S}} \subset \mathbb{S}}{\max } \quad \mu( \mathbb{\hat{S}}) + \gamma \sum_{i=1}^N \mathbb{I}(\hat{\xi}_i \in \hat{\mathbb{S}}) \\
~ & \text {s.t.} 
& f(x, \xi) \leq 0, \quad x \in \mathcal{X},	\quad \forall \xi \in \mathbb{\hat{S}}.
\end{array}
\end{equation}
Set $\mathbb{\hat{S}}$ to be the solution of stage $1$ gives stage $2$: 

\begin{equation}  \label{eq:stage_2_ori}
\begin{array}{lll}
\text{\emph{Stage 2}:}  & \underset{x}{\min } & g(x)  
\\
~ & \text {s.t.} 
& f(x, \xi) \leq 0, \quad x \in \mathcal{X},	\quad \forall \xi \in \mathbb{\hat{S}}.
\end{array}
\end{equation}
\begin{remark}
Note that when the new uncertainty set $\hat{\mathbb{S}}$ retains the same set representation as 
 $\mathbb{S}$,
 the problem \eqref{eq:stage_2_ori} can be reformulated in the same way for the original robust optimization problem \eqref{eq:problem_formulation} (see~\cite{beyer2007robust,ben2009robust,bertsimas2011theory} for more details). 
\end{remark}

In the following, we discuss some examples of uncertainty sets for which~\eqref{eq:stage_1} can be written as a finite program.
In particular, we scale the uncertainty set ${\mathbb{S}}$, and refer to the uncertainty sets introduced in~\cite{ben2015deriving} as follows.
\begin{itemize}
\item
\(p\)-norm  \(\|\cdot\|_p\),
$
\mathbb{S} = \{ \xi \mid \| \xi \|_p \leq 1 \}
$,
then
\eqref{eq:stage_1} becomes
\begin{equation} \label{eq:p_norm}
\begin{array}{ll}
\underset{x}{\max } & \underset{\alpha}{\max } \ \alpha +  \gamma \sum_{i=1}^N \mathbb{I}(\hat{\xi}_i \in \hat{\mathbb{S}})
\\
\text {s.t.} 
& f(x, \xi ) \leq 0, \quad x \in \mathcal{X}, 
\\& \forall \| \xi  \|_p \leq \alpha ,	\quad \alpha \in [0, 1].
\end{array}
\end{equation}
\item Polyhedron:
$
\mathbb{S} = \{ \xi  \mid V \xi  \leq \mathrm{1}\}
$
where \(V \in \mathbb{R}^{m \times n}\), \(1 \in \mathbb{R}^n\),
\eqref{eq:stage_1} becomes
\begin{equation} \label{eq:polyhedron}
\begin{array}{ll}
\underset{x}{\max } & \underset{v, \alpha }{\max } \ \alpha + \gamma \sum_{i=1}^N \mathbb{I}(\hat{\xi}_i \in \hat{\mathbb{S}})
\\
\text {s.t.} 
& f(x, \xi ) \leq 0, \quad x \in \mathcal{X},
\\& \forall \xi:  V \xi \leq \alpha \bm{1} + (1-\alpha) V v,
Vv  \leq 1, \alpha \in [0, 1]
\end{array}
\end{equation}
\item Cone-based set:
Let \( \mathbb{S} = \{ \xi  \mid \bm{1} - V \xi  \in C \} \), where \(V \in \mathbb{R}^{m \times n}\), \(\bm{1} \in \mathbb{R}^n\), and \( C \) is a pointed cone that contains a strictly feasible solution (i.e., there exists a \( \xi  \) such that \( \bm{1} - V \xi  \in \mathrm{int}(C) \)).
\eqref{eq:stage_1} becomes
\begin{equation} \label{eq:cone}
\begin{array}{ll}
\underset{x}{\max } & \underset{\alpha }{\max } \ \alpha + \gamma \sum_{i=1}^N \mathbb{I}(\hat{\xi}_i \in \hat{\mathbb{S}})
\\
\text {s.t.} 
& f(x, \xi ) \leq 0, \quad x \in \mathcal{X}, 
\\& \forall \xi:   \bm{1} - V \xi  \in  \alpha C, 
\ \alpha \in [0, 1]
\end{array}
\end{equation}
\item Semi-definite representable sets:
$\mathbb{S} = \{\xi \mid \lambda_{\max}(B(\xi )) \leq 1\}$, where \(B(\xi )\) is a matrix whose elements are linear in \(\xi \), and \(\lambda_{\max}\) denotes the maximum eigenvalue. This set can be reformulated as
$\mathbb{S} = \{\xi \mid \rho \mathrm{I} - B(\xi ) \succeq 0\}$, then
\eqref{eq:stage_1} becomes
\begin{equation} \label{eq:semi-definite}
\begin{array}{ll}
\underset{x}{\max } & \underset{\alpha }{\max } \ \alpha + \gamma \sum_{i=1}^N \mathbb{I}(\hat{\xi}_i \in \hat{\mathbb{S}})
\\
\text {s.t.} 
& f(x, \xi) \leq 0, \quad x \in \mathcal{X}, 
 \\& \forall \xi:  \alpha  \mathrm{I} - B(\xi) \succeq 0, \  \alpha \in [0, 1].
\end{array}
\end{equation}
\end{itemize}
\begin{remark}
If the robust counterpart $f(x, \cdot)$ is concave $\forall x \in \mathcal{X}$ and $f(\cdot, \xi)$ is convex $\forall  \xi \in \mathbb{S}$, then \eqref{eq:stage_2_ori} can be reformulated as a tractable optimization problem, for example the robust counterpart for a linear programming problem with polyhedral or ellipsoidal uncertainty regions can be reformulated as a linear programming or conic quadratic programming problem, we refer to~\cite{ben2009robust,ben2015deriving} for more details.
Besides, in \eqref{eq:p_norm}, \eqref{eq:cone} and \eqref{eq:semi-definite},  adding the indicator function will not influence the solution since larger uncertainty set, i.e., $\alpha$, includes more uncertainty samples, i.e., results in larger $\sum_{i=1}^N \mathbb{I}(\hat{\xi}_i \in \hat{\mathbb{S}})$. 
\end{remark}

\begin{figure}[h!]
    \centering
    \vspace{-1em}
    \includegraphics[width=0.475\linewidth]{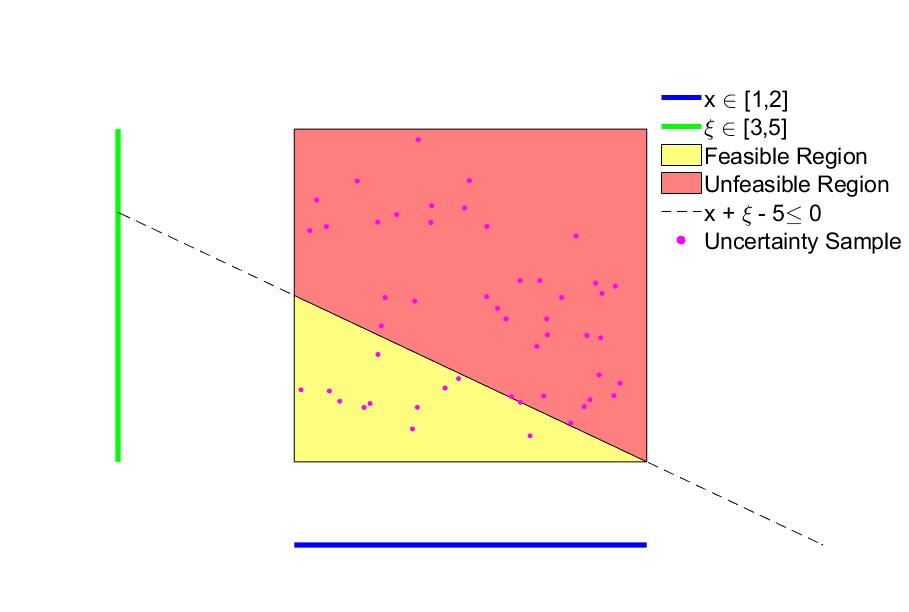}
    \includegraphics[width=0.475\linewidth]{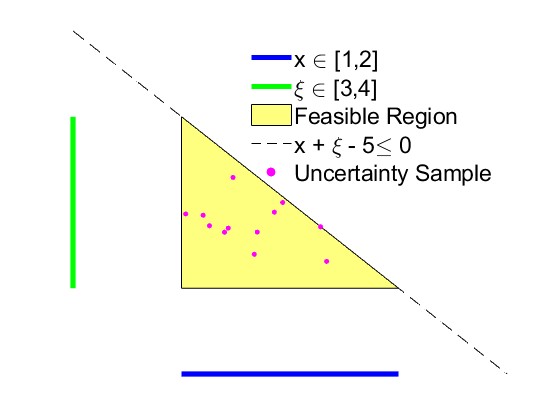}
    \vspace{-1em}
    \caption{Motivating example for \eqref{eq:stage_1}.
    On the left-hand side, we depict the problem with its original uncertainty set and the randomly generated uncertainty samples. On the right-hand side, we show the new uncertainty set that contains the most uncertainty samples while ensuring feasibility. 
    \vspace{-1em}
    }
    \label{fig:stage_1}
\end{figure}

\textbf{Motivating example revisited:} We use the same problem introduced in Section \ref{se:introduction} to illustrate~\eqref{eq:stage_1}, except that we randomly sample uncertain parameter $\xi$ from the set $\mathbb{S} = \{ \xi \mid  3 \leq \xi \leq 5 \}$, and the new uncertainty set obtained by solving \eqref{eq:stage_1} is still $\mathbb{\hat{S}} = \{ \xi \mid  3 \leq \xi \leq 4 \}$,
see Fig.~\ref{fig:stage_1}. \hfill $\blacktriangleleft$

\smallskip

\section{Data-driven distributionally adjustable robust optimization} \label{sec:unknown_uncertain}
In this section, we consider the case where the uncertainty $\xi$  is stochastic with an unknown distribution $\mathbb{P}$. 
We use the Wasserstein ambiguity set centered at the empirical distribution $\mathbb{P}^{\mathbb{N}}$ on $N$ i.i.d. training samples. 
Specifically, the Wasserstein distance of two distributions $\mathbb{Q}_1$ and $\mathbb{Q}_2$ 
represents the minimum cost required to transport the probability mass from $\mathbb{Q}_1$ to $\mathbb{Q}_2$. The corresponding Wasserstein ambiguity set includes all (continuous or discrete) distributions that are sufficiently close to the (discrete)  empirical distribution $\hat{\mathbb{P}}_N$ under this metric~\cite{mohajerin2018data}.
Given i.i.d. samples $\hat{\xi}^1, \dots, \hat{\xi}^N$ and the corresponding empirical distribution $\hat{\mathbb{P}}_N$, the  Wasserstein ball of radius $\epsilon$ is centered at the empirical distribution such that
$
\mathbb{B}_{\epsilon}(\mathbb{P}^{\mathbb{N}}) \coloneqq \{ \mathbb{Q} \in \mathcal{M}(\Xi) : d_{\text{Wass}}(\mathbb{P}, \hat{\mathbb{P}}_N) \leq \epsilon \}
$,
where $\mathcal{M}(\Xi)$  is the space of all probability distributions supported on $\Xi$,   $d_{\text{Wass}}(\mathbb{P}, \hat{\mathbb{P}}_N)$ denotes the Wasserstein distance between the distributions $\mathbb{P}$ and $\hat{\mathbb{P}}_N$. In stage $1$, to determine the best uncertainty set under which the constraints are satisfied with high probability,
we use the uncertainty quantification introduced in~\cite{mohajerin2018data}.
Let $\mathbb{Q}$ be the distribution of uncertainty,  the worst-case probability of the system being unsafe is 
\begin{equation} \label{eq:worst-case_init}
\underset{\mathbb{Q} \in \mathbb{B}_{\epsilon}(\mathbb{P}^{\mathbb{N}})}\max{\mathbb{Q}( \hat{\xi} \notin \hat{\mathbb{S}}} ),
\end{equation}
and the best-case probability of the system being safe is 
\begin{equation} \label{eq:best-case_init}
\underset{\mathbb{Q} \in \mathbb{B}_{\epsilon}(\mathbb{P}^{\mathbb{N}})}\max{\mathbb{Q}( \hat{\xi} \in \hat{\mathbb{S}} )},
\end{equation}
where $\epsilon$ is the Wasserstein distance, $ \mathbb{B}_{\epsilon}(\mathbb{P}^{\mathbb{N}})$ represents  the Wasserstein ball of radius  $\epsilon$ centered at the empirical distribution $\mathbb{P}^{\mathbb{N}}$~\cite{mohajerin2018data}.
By using the uncertainty quantifications \eqref{eq:worst-case_init} and \eqref{eq:best-case_init}, 
we now show that when the distribution is unknown,  we can reformulate the problem as finite programs as shown in \cite[Corollary 5.3]{mohajerin2018data}.

\begin{lemma} \label{lemma:safety}
Suppose that the uncertainty set is a polytope of the form $\mathbb{S} = \{\xi: C\xi \leq d\}$.
If the safety set $\mathbb{\hat{S}} = \{\xi: A\xi < b \in \mathbb{R}^L\}$ is an open polytope and the half-space $\lbrace \xi : \langle a_l, \xi \rangle \geq b_l \rbrace$ 
has a nonempty intersection with \( \mathbb{S} \) for any \( l \leq L \), where \( a_l \) is the \( l \)-th row of \( A \) and \( b_l \) is the \( l \)-th entry of the vector \( b \).
$\hat{\xi}_i, i=1,...,N$ are the uncertainty samples, $\epsilon > 0$ is the radius of the Wasserstein ball.
The worst-case probability \eqref{eq:worst-case_init} of the system being unsafe is given by
\begin{equation}  \label{eq:worst-case}
\begin{aligned}
\inf_{\lambda, s_i, \gamma_{il}, \theta_{il}} \quad & \lambda \epsilon + \frac{1}{N} \sum_{i=1}^{N} s_i \\
\operatorname{s.t.} \quad & 1 - \theta_{il} \left(b_l - \langle a_l, \hat{\xi}_i \rangle \right) + \langle \gamma_{il}, d - C \hat{\xi}_i \rangle \leq s_i, 
\\&
\forall i \leq N, \, \forall l \leq L \\
& \|{a}_l \theta_{il} - C^T \gamma_{il}\|_* \leq \lambda, \quad \forall i \leq N, \, \forall l \leq L \\
& \gamma_{i} \geq 0, \theta_{i} \geq 0, s_i \geq 0,   \quad \forall i \leq N, 
\end{aligned}
\end{equation}

If the safety set $\mathbb{\hat{S}} = \{\xi: A\xi \leq b \in \mathbb{R}^L\}$  is a closed polytope that has a nonempty intersection with $\mathbb{S}$, then the best-case probability of the system being safe is given by
\begin{equation} \label{eq:best-case}
\begin{aligned}
\inf_{\lambda, s_i, \gamma_i, \theta_i} \quad & \lambda \epsilon + \frac{1}{N} \sum_{i=1}^{N} s_i \\
\operatorname{s.t.} \quad & 1 + \langle \theta_i, \mathbf{b} - A \hat{\xi}_i \rangle + \langle \gamma_i, d - C \hat{\xi}_i \rangle \leq s_i \\& \forall i \leq N 
\\
& \|A^T \theta_i + C^T \gamma_i\|_* \leq \lambda, \quad \forall i \leq N \\
& \gamma_{i} \geq 0, \theta_{i} \geq 0, s_i \geq 0,   \quad \forall i \leq N.
\end{aligned}
\end{equation}
\end{lemma}

\smallskip
Now we can reformulate \eqref{eq:stage_1} as follows.
\subsubsection{Stage 1}
We choose the uncertainty set using the uncertainty quantification:
\begin{equation} \label{eq:uncertainty}
\begin{array}{ll}
 &\underset{x}{\min}   \underset{\mathbb{Q} \in \mathbb{B}_{\epsilon}(\mathbb{P}^{\mathbb{N}})}\max{\mathbb{Q}(\xi \notin \hat{\mathbb{S}} )} 
 \quad \ \text{or} \quad \ 
\underset{x}{\max}     \underset{\mathbb{Q} \in \mathbb{B}_{\epsilon}(\mathbb{P}^{\mathbb{N}})}\max{\mathbb{Q}(\xi \in \hat{\mathbb{S}})} 
\\[1em] 
& \text {s.t.} \quad
 \mathbb{\hat{S}}\subseteq \mathbb{S}, \quad f(x, \xi) \leq 0, 	\quad \forall \xi \in \mathbb{\hat{S}}.
\end{array}
\end{equation}

The following proposition presents the reformulation of the problem \eqref{eq:uncertainty}. 
\begin{proposition}
%\yulong{$L$ is not defined}
Assume that the uncertainty set is $\mathbb{S}=\{\xi  \in \mathbb{R}^{n}\mid  V \xi \leq \bm{1} \in \mathbb{R}^L\}$ and the assumptions of Lemma \ref{lemma:safety} holds. The worst case of being unsafe in \eqref{eq:uncertainty} can be rewritten in the  equivalent form \eqref{eq:worst-case-final}  
 \begin{equation} \label{eq:worst-case-final}
\begin{aligned}
&  \underset{\alpha, x, v, \lambda, s_i, \gamma_{i}, \theta_{i}}{\inf }     \quad  \lambda \epsilon + \frac{1}{N} \sum_{i=1}^{N} s_i
\\[-1em]&
\operatorname{s.t.}  \quad \quad
\begin{aligned}
 \\
 & 1 - \theta_{il} \left([\bm{1} + (1-\alpha) V v]_l - \langle V_l, \hat{\xi}_i \rangle \right) \\
 & \hspace{1cm }+ \langle \gamma_{il}, \bm{1} - V \hat{\xi}_i \rangle \leq s_i,  \quad  \forall i \leq N, \, \forall l \leq L \\
& \|{a}_l \theta_{il} - C^T \gamma_{il}\|_* \leq \lambda, \quad \forall i \leq N, \, \forall l \leq L \\
& \gamma_{i} \geq 0, \theta_{i} \geq 0,  s_i \geq 0, \quad \forall i \leq N,  
\\& 0\leq \alpha \leq 1,  V v \leq \bm{1}
\\& 
  f(x, \xi ) \leq 0, \quad x \in \mathcal{X},  \quad \forall \xi \in \mathbb{\hat{S}},
\end{aligned}
\end{aligned}
\end{equation}
where $L$ is defined as in Lemma~\ref{lemma:safety} and $\hat{\mathbb{S}}(\alpha, v) =\{ \xi \in \mathbb{R}^{n}\mid  V  \xi < \alpha \bm{1} + (1-\alpha) V v, \alpha \in [0,1], v \in \mathbb{S}\}$.
The best case of being safe in \eqref{eq:uncertainty} can be equivalently rewritten as 
\begin{equation} \label{eq:best-case-final}
\begin{aligned}
&  \underset{\alpha, x, v, \lambda, s_i, \gamma_{i}, \theta_{i}}{\max } \quad z 
\\[-2em]&
\operatorname{s.t.}  \quad \quad
\begin{aligned}
\\& 
z \leq \lambda \epsilon + \frac{1}{N} \sum_{i=1}^{N} s_i
 \\
 & 1 + \langle \theta_i, \alpha \bm{1} + (1-\alpha) V v - V \hat{\xi}_i \rangle \\
 & \hspace{1cm }+ \langle \gamma_i, \bm{1} - V \hat{\xi}_i \rangle \leq s_i, \quad
 \forall i \leq N \\
& \|A^T \theta_i + V^T \gamma_i\|_* \leq \lambda, \quad \forall i \leq N \\
& \gamma_i \geq 0, \theta_i \geq 0, s_i \geq 0,\quad \forall i \leq N 
\\& 0\leq \alpha \leq 1,  V v \leq \bm{1}
\\& 
 f(x, \xi ) \leq 0, \quad x \in \mathcal{X},  \quad \forall \xi \in \mathbb{\hat{S}}.
\end{aligned}
\end{aligned}
\end{equation}
\end{proposition}
\begin{proof}
We first consider the inner problem of \eqref{eq:uncertainty}.
Define the safety set
$\hat{\mathbb{S}}(\alpha, v) =\{ \xi \in \mathbb{R}^{n}\mid  V  \xi < \alpha \bm{1} + (1-\alpha) V v, \alpha \in [0,1], v \in \mathbb{S}\}$, note that $\hat{\mathbb{S}}(\alpha, v)\subseteq \mathbb{S}$ for all $\alpha \in [0,1]$ and $v\in \mathbb{S}$. Replacing $\Xi $ and $\mathbb{A}$ with $\mathbb{S}$ and $\hat{\mathbb{S}}(\alpha, v) $, respectively, it follows that the assumptions in Lemma \ref{lemma:safety} hold. According to Lemma \ref{lemma:safety}, the worst-case \eqref{eq:worst-case} is written as 
\begin{equation}  \label{eq:worst-case-inf}
\begin{aligned}
\inf_{\lambda, s_i, \gamma_{i}, \theta_{i}} \quad & \lambda \epsilon + \frac{1}{N} \sum_{i=1}^{N} s_i \\
\text{s.t.} \quad & 1 - \theta_{il} \left([\alpha \bm{1} + (1-\alpha) V v]_l - \langle V_l, \hat{\xi}_i \rangle \right) 
\\& + \langle \gamma_{il}, 1 - V \hat{\xi}_i \rangle \leq s_i, \quad \forall i \leq N, \, \forall l \leq L \\
& \|\mathbf{a}_l \theta_{il} - C^T \gamma_{il}\|_* \leq \lambda, \quad \forall i \leq N, \, \forall l \leq L \\
& \gamma_{i} \geq 0, \theta_{i} \geq 0, s_i \geq 0,   \quad \forall i \leq N.
\end{aligned}
\end{equation}
where $[\bm{1} + (1-\alpha) V v]_l$ and $V_l$ are the $l-$th entry and $l-$th row of the vector $\bm{1} + (1-\alpha) V v$ and the matrix $V$.

Similarly, 
the best-case \eqref{eq:best-case} is written as
\begin{equation} \label{eq:best-case-inf}
\begin{aligned}
\inf_{\lambda, s_i, \gamma_i, \theta_i} \quad & \lambda \epsilon + \frac{1}{N} \sum_{i=1}^{N} s_i \\
\text{s.t.} \quad & 1 + \langle \theta_i, \alpha \bm{1} + (1-\alpha) V v - V \hat{\xi}_i \rangle \\
& \hspace{1cm }+ \langle \gamma_i, \bm{1} - V \hat{\xi}_i \rangle \leq s_i, \quad \forall i \leq N \\
& \|A^T \theta_i + V^T \gamma_i\|_* \leq \lambda, \quad \forall i \leq N \\
& \gamma_{i} \geq 0, \theta_{i} \geq 0, s_i \geq 0,   \quad \forall i \leq N.
\end{aligned}
\end{equation}

Therefore, when using the worst-case probability,
\eqref{eq:uncertainty} is rewritten as  \eqref{eq:worst-case-final}.
When using the best-case probability,
\eqref{eq:uncertainty} is rewritten as
\begin{equation}
\begin{aligned}
&  \underset{\alpha, x, v}{\max }    \inf_{\lambda, s_i, \gamma_{i}, \theta_{i}} \quad  \lambda \epsilon + \frac{1}{N} \sum_{i=1}^{N} s_i
\\[-1em]&
\text {s.t.}  \quad \quad
\begin{aligned}
 \\
 & 1 + \langle \theta_i, \alpha \bm{1} + (1-\alpha) V v - V \hat{\xi}_i \rangle \\
 & \hspace{1cm}+ \langle \gamma_i, \bm{1} - V \hat{\xi}_i \rangle \leq s_i, \quad \forall i \leq N \\
& \|A^T \theta_i + V^T \gamma_i\|_* \leq \lambda, \quad \forall i \leq N \\
& \gamma_{i} \geq 0, \theta_{i} \geq 0, s_i \geq 0,   \quad \forall i \leq N
\\& 0\leq \alpha \leq 1,  V v \leq \bm{1}
\\& 
  f(x, \xi ) \leq 0, \quad \forall \xi \in \mathbb{\hat{S}}.
\end{aligned}
\end{aligned}
\end{equation}
which can be written as \eqref{eq:best-case-final}.
\end{proof}

\subsubsection{Stage 2}
Fix the uncertainty set $\hat{\mathbb{S}}(\alpha, v) =\{ \xi \in \mathbb{R}^{n}\mid  V  \xi < \alpha^* \bm{1} + (1-\alpha^*) V v^*\} $ based on   the solution $\alpha^*$, $v^*$ of stage $1$, solve the robust optimization problem~\eqref{eq:stage_2_ori}.

\section{Numerical experiments}
In this section, we consider a numerical example and an example of an optimal power flow problem. Problems \eqref{numerical:eq_stege_1}, \eqref{numerical:distributionally_robust}, \eqref{simulation:eq_stege_1} and \eqref{simulation:distributionally_robust} are solved by \text{IPOPT}~\cite{wachter2006implementation}.

\subsection{Numerical example}
Consider the following optimization problem
\begin{equation} \label{eq:numerical_problem}
\begin{aligned}
\min_{x \in \mathbb{R}^2} & \quad \bm{c}^T \bm{x}
\\
\operatorname{s.t.} & \quad A \bm{x} \leq \bm{b} + \xi
\\&
\quad \bm{x} \geq 0, \quad \forall B\xi \leq  \bm{d}
\end{aligned}
\end{equation}
where $x$ is the decision variable, $\xi$ is the uncertain parameter, $c = \begin{bmatrix}
    -1;-3
\end{bmatrix}$
$A =\begin{bmatrix}
    1 & 1 \\1 & 2
\end{bmatrix}$, $ \bm{b} = \begin{bmatrix}
    5\\6
\end{bmatrix}$, $B =\begin{bmatrix}
    1 & 0; 0 & 1;-1 & 0; 0& -1
\end{bmatrix}$, $\bm{d} = \begin{bmatrix}
    10;10;10;10
\end{bmatrix}$.

\begin{comment}
\begin{figure}[t!]
\centering
\vspace{-1em}
\includegraphics[width=1\linewidth]{case_feasible_region.png}
\vspace{-3em}
\caption{Feasible region for the optimization problem~\eqref{eq:numerical_problem}.}
\label{fig:case_feasible_region}
\vspace{-1em}
\end{figure}
\end{comment}
\subsubsection{The case of non-stochastic uncertainty} When the uncertainty is non-stochastic, denote $\mathbb{\hat{S}} = 
\{ \xi  \in \mathbb{R}^{2}\mid  B \xi \leq \alpha \bm{d} + (1-\alpha) B v,  B v \leq \bm{d}
\}$, we refer to Section \ref{sec:uniform_uncertain} and rewrite the stage $1$ of the original problem \eqref{eq:numerical_problem} as follows
\begin{equation}  \label{numerical:eq_stege_1}
\begin{aligned}
& \underset{v, \alpha, x}{\max } \ \alpha + \gamma \sum_{i=1}^N \mathbb{I}(\hat{\xi}_i \in \hat{\mathbb{S}})
\\[-1em]&
\text {s.t.}  \quad \quad
\begin{aligned} 
\\&  A \bm{x} \leq \bm{b} + \xi, 
\quad  x \geq 0
\\& \forall  B\xi  \leq \alpha \bm{d} + (1-\alpha) B v
\\&\alpha \in [0,1], \quad  B v \leq \bm{d},
\end{aligned}
\end{aligned}
\end{equation}
where the $500$ uncertainty samples $\hat{\xi}_{i=1,...,500}$ are randomly generated from $\mathbb{{S}} = \{ \xi  \in \mathbb{R}^2 \mid  B \xi  \leq   \bm{d}  \}$.
The solution of \eqref{numerical:eq_stege_1} is $\alpha^* = \frac{3}{4}$, $v^* = 10$, the new uncertainty set is $\mathbb{\hat{S}} = 
\{ \xi  \in \mathbb{R}^{2}\mid  B \xi \leq  \begin{bmatrix}
    10;10;5;5
\end{bmatrix}
\}$, and the solution of \eqref{eq:numerical_problem} with new uncertainty set is $x = \begin{bmatrix}
    0;0
\end{bmatrix}$.
\smallskip

\subsubsection{The case of stochastic uncertainty} When the uncertainty set is stochastic but with an unknown distribution, we refer to Section \ref{sec:unknown_uncertain} and 
the safety set would be $ \hat{\mathbb{S}}$. 
Using the best-case probability, stage $1$ of the original problem \eqref{eq:numerical_problem} is written as follows
\begin{equation} \label{numerical:distributionally_robust}
\begin{aligned}
&  \underset{\alpha, v, x, \lambda, s, \gamma, \theta}{\max } \quad z 
\\&
\text {s.t.}  \quad
\begin{aligned}
\\[-2em]& 
z \leq \lambda \epsilon + \frac{1}{N} \sum_{i=1}^{N} s_i
 \\
 & 1 + \langle \theta_i, \alpha \bm{d} + (1-\alpha) B v - B \hat{\xi}_i \rangle + \langle \gamma_i, \bm{d} - B \hat{\xi}_i \rangle \leq s_i, 
 \\&
 \forall i \leq N 
 \\&  \| B^T \theta_i + B^T\gamma_i\| \leq \lambda, \quad \forall i \leq N 
\\& \gamma_i \geq 0, \theta_i \geq 0, s_i \geq 0,\quad \forall i \leq N \\ 
& \forall  B\xi  \leq \alpha \bm{d} + (1-\alpha) B v,
\\&  A \bm{x} \leq \bm{b} + \xi, 
\quad  x \geq 0, \quad   \alpha \in [0,1], \quad  B v \leq \bm{d}, \\
\end{aligned}
\end{aligned}
\end{equation}
We select the Wasserstein radius  $\epsilon = 0.05$ and   $N = 100$, the new uncertainty set becomes $\hat{\mathbb{S}} = \{ \xi  \in \mathbb{R}^{2}\mid  B \xi \leq  \begin{bmatrix}
    8.1;7.7;5.0;5.5
\end{bmatrix}
\}$, when $N=500$, the new uncertainty set is the same with case 1), i.e.,   $\hat{\mathbb{S}} = \{ \xi  \in \mathbb{R}^{2}\mid  B \xi \leq  \begin{bmatrix}
    10;10;5;5
\end{bmatrix}
\}$.

\begin{figure}[t]
    \centering
\includegraphics[width=0.6\linewidth]{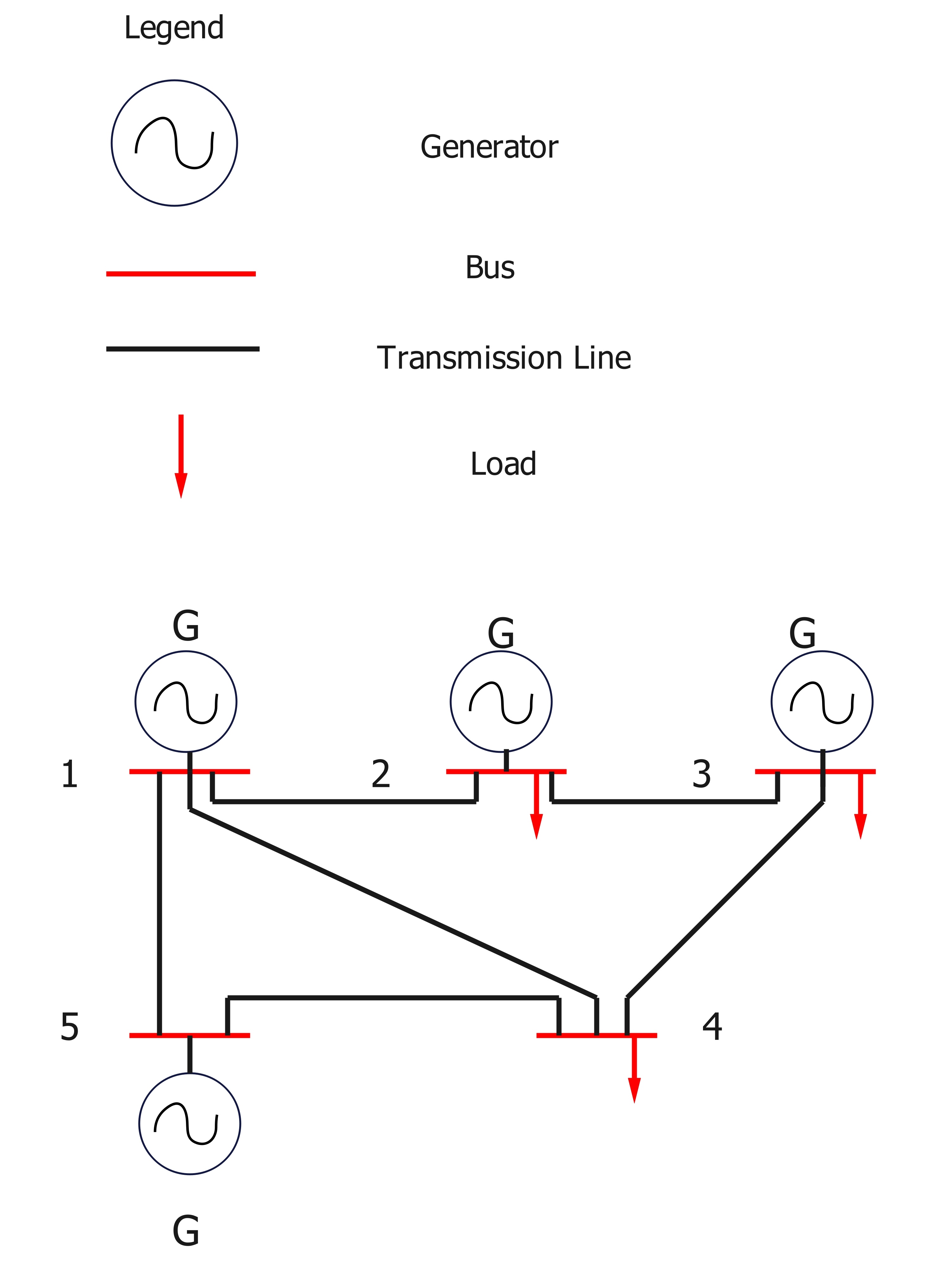}
\vspace{-1em}
\caption{IEEE-5-bus-system}
    \label{fig:ieee_bus}
\vspace{-1em}
\end{figure}
\subsection{Case study in optimal power flow}
In this section, we illustrate the theoretical results of this paper with an optimal power flow problem under uncertainty.
Case studies are carried out through the modified IEEE 5-bus distribution system~\cite{PowerSystemsTestCaseArchive} as shown in Fig.~\ref{fig:ieee_bus}.
Some system parameters are set as follows: total load demand, \(\sum_{d_i\in\mathcal{D}} P_{d_i} = 1000\text{MW}\), 
the largest loss of generation, $P_L = \max_{g_i\in\mathcal{G}} P_{g_i} $, where $\mathcal{D}$ and $\mathcal{G}$ are the sets of loads and generators respectively, $P_{g_1}$ and $P_{g_2}$ are renewable energy generators; Base power, $S_{\text{base}} =  100 \text{MVA}$, PFR delivery time \( T_d = 10 \, \text{s} \), total PFR is $ R =  324\text{MW}$. 
The frequency limits of nadir and the rate-of-change-of-frequency are set as: \( \Delta f_{\text{lim}} = 0.8 \, \text{Hz} \), and \( \Delta \dot{f}_{\text{lim}} = 0.5 \, \text{Hz/s} \). The power system network is defined by all the nodes ($n\in\mathcal{N}$) and branches ($mn \in\mathcal{R}$).
Other parameters are available at~\cite{MATPOWERDocumentation}.

In what follows, the optimal power flow problem is formulated with the aim of minimizing the total generation cost and the amount of virtual inertia ($H_{c_i}$) provided by the inverter-based resources. The constraints of frequency nadir and the rate-of-change-of-frequency are defined in \eqref{eq:nadir} and \eqref{eq:rocof} respectively, with the inertia $H_i$ being the uncertain parameter. Power balance constraints at each node in the system are given by \eqref{eq:power_balance} and power flow within each line is defined in \eqref{eq:power_flow}. The remaining constraints represent the bounds of the decision variables and uncertain parameters.
\begin{subequations} \label{p:opf}
 \begin{align}
\underset{H_c, \theta, P_{g_{i}}}{\min} \quad & J = \sum_{g_i\in\mathcal{G}} c_{g_i} P_{g_i}  + \sum_{c_i\in\mathcal{C}} c_{c_i} H_{c_i}
\\ \text{s.t.}:  \quad  &  HR\ge \frac{\Delta P_L^2T_d}{4\Delta f_\mathrm{lim}} \label{eq:nadir} \\&    \frac{P_L}{2H}  \leq \Delta \dot{f}_{\text{lim}}  \label{eq:rocof}
\\& H = \sum_{g_i\in\mathcal{G}} H_{g_i} + \sum_{c_i\in\mathcal{C}} H_{c_i} + H_{i}
\\&  \sum_{g_i\in\mathcal{G}_n}P_{g_i,n} - \sum_{d_i\in\mathcal{D}_n}P_{d_i,n} = \sum_{mn\in\mathcal{R}} P_{mn},\, \forall n\in\mathcal{N} \label{eq:power_balance}
\\&
P_{mn} =  \frac{\theta_m-\theta_n}{X_{mn}},\, \forall mn\in\mathcal{R} \label{eq:power_flow}
\\& 
\underline {P}_{g_i} \le P_{g_i} \leq \Bar{P}_{g_i} ,\, \forall g_i\in\mathcal{G} \\&
0\leq \theta_n \leq \pi,\, \forall n\in\mathcal{N}
\\& \underline {H}_{c} \leq H_{c_i} \leq \Bar {H}_{c},\, \forall C_i\in\mathcal{C}, \underline {H}_{i} \leq H_i \leq \Bar {H}_{i}
\end{align}   
\end{subequations}

In the following, we summarize the constraint in \eqref{p:opf} as $ f(x, H_i ) \leq 0, \quad x \in \mathcal{X},$ in which $x$ represents the decision variables and $H_i$ is the uncertain parameter.
Denote $A = \begin{bmatrix}
1 ; -1
\end{bmatrix}$,  $b = \begin{bmatrix}
35 ; -20
\end{bmatrix}$.

\smallskip

\subsubsection{The case of non-stochastic  uncertainty}
When the uncertainty is non-stochastic, denote $\mathbb{\hat{S}} = 
\{ H_i  \in \mathbb{R}\mid  A  H_i \leq \alpha \bm{b} + (1-\alpha) A v,  A v \leq \bm{b}
\}$, $ H_{c} = \sum_{c_i\in\mathcal{C}} H_{c_i}$, we refer to Section \ref{sec:uniform_uncertain} and rewrite the stage $1$ of the original problem \eqref{p:opf} as follows
\begin{equation}  \label{simulation:eq_stege_1}
\begin{aligned}
& \underset{v, \alpha, H_c}{\max } \ \alpha + \gamma \sum_{i=1}^N \mathbb{I}(\hat{\xi}_i \in \hat{\mathbb{S}})
\\[-2em]&
\text {s.t.}  \quad \quad
\begin{aligned} 
\\&  (600 + H_c+ H_i) 324 \geq \frac{829440}{3.2} 
\\&  \frac{600}{2(600 + H_c+ H_i)}  \leq 0.5  
\\& 116 \leq H_c \leq 175
\\&  A H_i \leq \alpha \bm{b} + (1-\alpha) A v, \quad  A v \leq b, \quad \alpha \in [0,1]
\end{aligned}
\end{aligned}
\end{equation}
where the $500$ uncertainty samples $\hat{\xi}_{i=1,...,500}$ are randomly generated from $\mathbb{{S}} = \{ H_i  \in \mathbb{R}\mid 20 \leq H_i \leq 35\}$.

The solution of \eqref{simulation:eq_stege_1} is $\alpha^* = \frac{2}{3}$, $v^* = 35$. So the new uncertainty set is 
$
\mathbb{\hat{S}} = 
\{ H_i  \in \mathbb{R}\mid 25 \leq H_i \leq 35
\}$, and it is easy to verify that it is the largest subset that makes the problem feasible. In stage $2$, we use the new uncertainty set of $H_i$ and solve \eqref{p:opf}. 
The solution of \eqref{p:opf} with new uncertainty set $\mathbb{\hat{S}}$ is $P_{g_1} = 40, P_{g_2} = 170, P_{g_3} = 520, P_{g_5} = 270,  H_{c_1} = 100, H_{c_2} = 75, \theta_1 = 0.0373, \theta_2 = 0.0272, \theta_3 = 0.0374, \theta_4 = 0, \theta_5 = 0.0449$.

\smallskip

\subsubsection{The case of stochastic uncertainty}
When the uncertainty set is stochastic but with an unknown distribution, we refer to Section \ref{sec:unknown_uncertain} and
the safety set would be $\hat{\mathbb{S}}$. 
By using the best-case probability, the problem \eqref{p:opf} can be written as 
\begin{equation*} \label{simulation:distributionally_robust}
\begin{aligned}
&  \underset{\alpha, v, H_c , \lambda, s, \gamma, \theta}{\max } \quad z 
\\[-2em]&
\text {s.t.}  \quad \quad
\begin{aligned}
\\& 
z \leq \lambda \epsilon + \frac{1}{N} \sum_{i=1}^{N} s_i
 \\
 & 1 + \langle \theta_i, \alpha \bm{b} + (1-\alpha) A v - A \hat{\xi}_i \rangle + \langle \gamma_i, \bm{b} - A \hat{\xi}_i \rangle \leq s_i, 
 \\&
 \forall i \leq N 
 \\
& \| A^T \theta_i + A^T\gamma_i\| \leq \lambda, \quad \forall i \leq N \\
& \gamma_i \geq 0, \theta_i \geq 0, s_i \geq 0,\quad \forall i \leq N 
\\&  (550 + H_c+ H_i) 300 \geq \frac{720000}{3.2} 
\\&  \frac{600}{2(550 + H_c+ H_i)}  \leq 0.5  
\\& 116 \leq H_c \leq 175
\\& A H_i \leq \alpha \bm{b} + (1-\alpha) A v,  \quad A v \leq b, \quad \alpha \in [0,1].
\end{aligned}
\end{aligned}
\end{equation*}
We select the Wasserstein radius from $\epsilon  \in \{0.1, 0.05\}$, the cardinality of training datasets from  $N \in \{100, 500, 1000\}$, the new uncertainty set under different parameters are provided in Table \ref{table:case_2}, it is shown that in this case, choosing the radius $\epsilon = 0.05$ provides larger uncertainty set.

\begin{table}[t!]
\centering
\scalebox{1.05}{ 
\begin{tabular}{|c|c|c|c|c|}
\hline
\multicolumn{2}{|c|}{\raisebox{0\height}{$\mathbb{\hat{S}}$}} & \multicolumn{3}{c|}{$N$} \\
\cline{3-5}
\multicolumn{2}{|c|}{} & 100 & 500 & 1000 \\
\hline
\multirow{2}{*}{$\epsilon$} & $0.1$ &  [25.56,\ 30.77] & [25.89,\ 34.67]   & [25.00,\ 34.6]\\
\cline{2-5}
&  $0.05$ & [25.01,\ 34.78] & [25.00,\ 34.79] &[25.11,\ 34.99] \\
\hline
\end{tabular}
    }
\caption{New uncertainty set under different parameters $N$ and $\epsilon$.}
\vspace{-3em}

\label{table:case_2}
\end{table}

\smallskip

\section{Concluding Remarks}
This paper investigates robust optimization with adjustable uncertainty sets. In particular, we have proposed to utilize data-driven methods to determine the largest uncertainty set under which the problem is feasible and, then, solve the robust optimization problem with the devised uncertainty set. 
For non-stochastic uncertainty, we have identified a scaled uncertainty set that includes the maximum number of uncertainty samples, demonstrating that this reformulation is straightforward for commonly used uncertainty sets. For stochastic uncertainty, we have performed uncertainty quantification over the ambiguity sets defined by Wasserstein metric. In turn, when the uncertainty is a polytope, we have further shown that the problem can be reformulated as a finite program.
The effectiveness of the proposed approaches has been validated through numerical examples and a case study on an optimal power flow problem under uncertainty.

%\balance
\bibliographystyle{ieeetr}

\end{document}